\documentclass[a4paper]{amsart}
\usepackage[all]{xy}
\usepackage{hyperref}
\usepackage{enumerate}

\usepackage{amsmath}
\usepackage{amsthm}
\usepackage{amssymb}
\usepackage{amscd}
\usepackage{graphicx}
\usepackage{epsfig}
\usepackage[all]{xy}

\numberwithin{equation}{section}
\allowdisplaybreaks

\newtheorem{theorem}{Theorem}[section]

\newtheorem{proposition}{Proposition}[section]
\newtheorem{corollary}{Corollary}[section]

\newtheorem{conjecture}{Conjecture}[section]
\theoremstyle{definition}
\newtheorem{definition}{Definition}[section]
\newtheorem{example}{Example} [section]

\newtheorem{remark}{Remark}[section]

\newcommand{\Rr}{\mathbb R}
\newcommand{\Ss}{\mathbb S}

\newcommand{\Nn}{\mathbb N}
\newcommand{\Tt}{\mathbb T}
\newcommand{\Qq}{\mathbb Q}
\renewcommand{\d}{\mathrm d}

\newcommand{\set}[1]{\left\{#1\right\}}

\newcommand{\X}{\ensuremath{\mathfrak{X}}}
\newcommand{\F}{\ensuremath{\mathcal{F}}}

\newcommand{\Po}{\text{Poiss}}

\newcommand{\R}{\mathcal{R}}

\newcommand{\G}{\mathcal{G}}            
\newcommand{\M}{\mathcal{M}}          
\newcommand{\tto}{\rightrightarrows}    
\newcommand{\Diff}{\text{\rm Diff}\,}   

\newcommand{\al}{\alpha}                
\newcommand{\be}{\beta}                 
\newcommand{\Lie}{\mathcal{L}}          
\renewcommand{\gg}{\mathfrak{g}}        

\begin{document}
\title{Rigidity and Flexibility in Poisson Geometry}

\author{Marius Crainic}
\address{Depart. of Math., Utrecht University, 3508 TA Utrecht, 
The Netherlands}
\email{crainic@math.uu.nl}
\thanks{Supported in part by NWO and a Miller Research Fellowship}

\author{Rui Loja Fernandes}
\address{Depart.~de Matem\'{a}tica, 
Instituto Superior T\'{e}cnico, 1049-001 Lisboa, PORTUGAL} 
\email{rfern@math.ist.utl.pt}
\thanks{Supported in part by FCT/POCTI/FEDER and by grant
  POCTI/1999/MAT/33081.}

\begin{abstract}
We study rigidity and flexibility phenomena in the context of
Poisson geometry.
\end{abstract}

\maketitle

\section{Introduction}
Rigidity and flexibility phenomena are well-known in symplectic
geometry. Since a Poisson manifold is a generalization of a symplectic
manifold, it is natural to look for rigidity and flexibility phenomena
in the context of Poisson geometry. In the present work we will survey
recent results in this direction, give a few new results, and pose some
conjectures. It should came as no surprise that these kind of questions are
even harder to investigate in this more general context.

In symplectic geometry there are two basic results, going in opposite
directions:
\begin{itemize}
\item For \emph{open manifolds}, a celebrated theorem by Gromov (\cite{Gro})
  states that every non-degenerate two form is homotopic to a symplectic
  form (\emph{flexibility}). 
\item For \emph{compact manifolds}, Moser's stability theorem (see,
  e.g., \cite{MaSa} Theorem 3.17) shows that two symplectic forms
  that are homotopic in the same cohomology class must be
  symplectomorphic (\emph{rigidity}).
\end{itemize}
These two results have important consequences. For example, the first
result reduces the problem of existence of a symplectic form on an open
manifold to the problem of existence of a non-degenerate 2-form, which
is a relatively easy problem in obstruction theory. In contrast, for
compact 4-manifolds there are highly non-trivial obstructions to the
existence of symplectic forms. In Poisson geometry, one would like to
have generalizations of these results, and in this paper we discuss
some recent progress in this direction.

In Section \ref{sec:flex}, we start by recalling Bertelson's theorem
(\cite{Ber}), which gives a partial analogue of Gromov's result in
Poisson geometry. It states that on an \emph{open foliated manifold} every
leafwise non-degenerate 2-form is homotopic, in the class of
non-degenerate 2-forms, to a leafwise symplectic form, i.e., a Poisson
structure. Therefore, for an \emph{open foliation} $\F$ the problem of
existence of a Poisson structure whose underlying foliation is $\F$,
can be reduced to a problem in obstruction theory. We point out that
Bertelson's results can be better understood in the setting of
groupoids. Unfortunately, these results are only valid for regular
foliations, while in Poisson geometry, some of the most interesting
Poisson manifolds have singular symplectic foliations. It would be
very interesting to find generalizations of these results to singular
foliations, but no such result is available at present.

In Section \ref{sec:rigid}, we turn to rigidity. Given a Poisson
structure $\pi$, one may ask if every nearby Poisson structure,
satisfying some cohomological condition, must be diffeomorphic to
$\pi$. We don't know of such analogue of Moser's stability theorem in
Poisson geometry. In fact, it is not clear how such a general result
should be stated since there is no obvious way of comparing Poisson
cohomology classes attached to distinct Poisson structures. Even the
question of what is the right compactness assumption is an interesting
problem. We will see that one is led naturally to the study of
\emph{proper} Poisson manifolds (Poisson manifolds which integrate to
proper symplectic groupoids) and of Poisson manifolds of \emph{compact
type} (Poisson manifolds which integrate to compact symplectic
groupoids). Properness is enough to guarantee rigidity around leafs,
as illustrated by recent results of Weistein and Zung on linearization
around leaves (\cite{Wein2,Zun1}). Some kind of global rigidity should
also hold for Poisson manifolds of compact type. However this is an
intriguing class of which very little seems to be known. We will give
some properties of this class, and we will state a rigidity result in
the regular case.

In face of such difficulties, one is led naturally to the study of
weaker properties of Poisson manifolds, but still keeping a flavor of
rigidity/flexibility. In this direction, one very natural question
which does not seem to have been discussed before, is the problem of
\emph{stability} of symplectic leaves: a leaf of a Poisson structure
is called stable if every nearby Poisson structure has a nearby
diffeomorphic leaf. In Section \ref{sec:stab}, we discuss the
following result which gives a sufficient condition for stability for
fixed points, i.e., zero dimensional symplectic leafs:

\begin{theorem}
Let $x_0\in M$ be a fixed point of a Poisson structure $\pi$ on $M$. 
If the isotropy Lie algebra $\gg_{x_0}$ has vanishing second
Chevalley-Eilenberg cohomology, i.e., $H^2(\gg_{x_0})=0$, then $x_0$
is a stable fixed point. 
\end{theorem}

We will also argue that this is the best possible general result in
this direction, and we will formulate a conjecture concerning
stability of symplectic leafs of higher dimension, which uses the
relative Poisson cohomology of the leaf due to Ginzburg and Lu
\cite{GiLu}.

\section{Flexibility}
\label{sec:flex}

In this section we start by recalling the results due to Melanie Bertelson
\cite{Ber}, which extend Gromov's results to the foliated case.

Recall that a manifold $M$ is open if and only if it carries a
positive, proper Morse function, without any local maximum. The
existence of such a Morse function is precisely what is need in order
to proof Gromov's h-principle for open invariant relations (see, e.g.,
\cite{Gei}, Theorem 3.3). Loosely speaking, the h-principle states
that for every open, invariant relation $\R$ on a jet space $J^kE$,
the space of sections $\Gamma(\R)$ is weak homotopy equivalent to the
space of holonomic sections $\Gamma_\text{hol}(\R)$.  The h-principle,
in turn, when applied to the relation
\[ \R=\set{j^1\al(x)\in J^1T^*M:\d\al(x) \text{ is non-degenerate}},\]
yields the following important result: 

\begin{theorem}[Gromov \cite{Gro}]
Every non-degenerate two form is homotopic in the set of
non-degenerate two forms to a symplectic form. 
\end{theorem}

Let now $\F$ be a foliation of a manifold $M$, and denote by
$\Omega^2(\F)$ the space of foliated 2-forms on $\F$. Note that a
leafwise symplectic 2-form $\omega\in\Omega^2(\F)$ is the same thing
as a Poisson structure on $M$ with symplectic foliation equal to
$\F$. In order to obtain an analogue of Gromov's result, one needs
the concept of an open foliation:

\begin{definition}[Bertelson \cite{Ber}]
\label{defn:open:foliation}
A foliated manifold $(M,\F)$ is said to be \textbf{open} if there exists
a smooth function $f: M \to [0,\infty)$with the following properties: 
\begin{enumerate}
\item[(a)] $f$ is proper;
\item[(b)] $f$ has no leafwise local maxima; 
\item[(c)] $f$ is $\F$-generic.
\end{enumerate}
\end{definition}

The condition that $f$ is $\F$-generic is a transversality condition
on its jet which plays a role analogous to the Morse condition.  The
main result of \cite{Ber} then states that on an open foliated
manifold, any open, foliated, invariant differential relation
satisfies the parametric h-principle.  Applying this result to the
relation
\[ \R=\set{j^1\al(x)\in J^1T^*\F:\d_\F\al(x) \text{ is non-degenerate}},\]
one obtains a foliated version of Gromov's theorem:

\begin{theorem}[Bertelson \cite{Ber}]
Let $(M,\F)$ be an open foliated manifold.  Any leafwise nondegenerate
$2$-form is homotopic, in the class of leafwise nondegenerate
$2$-forms, to a leafwise symplectic form.
\end{theorem}

The leaves of an open foliation are necessarily open manifolds. However,
this condition is not sufficient for the conclusions of the theorem
to hold. The following example is presented in \cite{Ber2}.

\begin{example}
Let $M$ be a compact $m$-manifold and $\F$ a dimension $2k$ foliation
of $M$.  If $\F$ is transversely orientable then $\F$ does not carry
an exact leafwise symplectic form. In fact, let $\mu$ be a transverse
volume form, i.e., a nowhere vanishing, closed, $(m-2k)$-form such
that $i_X\mu=0$, for any tangential vector field $X\in\X(\F)$. If
$\omega=\d_\F\al$ is an exact leafwise symplectic form, we choose
extensions $\widetilde{\omega}$, $\widetilde{\al}$ of $\omega$ and
$\al$ to forms on $M$ satisfying
$\widetilde{\omega}=\d\widetilde{\al}$, and we observe that:
\[\mu\wedge\widetilde{\omega}^k
=\d(\mu\wedge\widetilde{\omega}^{k-1}\wedge\widetilde{\al})\]
is an exact volume form on $M$. This contradicts the assumption that
$M$ was compact.

{For} example, take $M=\Tt^2\times\Ss^3$ with the foliation
$\F=\F_a\times\Ss^3$, where $\F_a$ is the irrational foliation of
the torus of slope $a$. By this we mean the foliation of $\Tt^2$
generated by the vector field $X_a=\frac{\partial}{\partial\theta_1}+
a\frac{\partial}{\partial\theta_2}$, where $(\theta_1,\theta_2)$ are
coordinates on $\Tt^2$ and $a\in\Rr-\Qq$. Now observe that:
\begin{enumerate}
\item[(a)] $\F$ carries a leafwise non-degenerate 2-form: let $\al$
  be the 1-form on $\Tt^2$ such that $\al(X)=1$, and let $\beta$
  be the contact form on $\Ss^3$. Then $\al\wedge\be+\d\be$ gives a
  leafwise non-degenerate 2-form on $\F$.
\item[(b)] $\F$ is transversely orientable: a transverse volume form
  is given by the 1-form $\mu=\pi_1^*(a\d\theta_1-\d\theta_2)$, where
  $\pi_1:\Tt^2\times\Ss^3\to\Tt^2$ denotes projection to the first
  factor.
\end{enumerate}
Therefore, $\F$ has open leafs and admits a leafwise non-degenerate
2-form. However, it does not carry any exact leafwise symplectic form.
\end{example}
Several other examples are given by M.~Bertelson in \cite{Ber2},
including examples with non-compact $M$. It should be observed that
in all such examples one has $H^2(\F)=0$. For instance, in the example
above this follows from:
\[ H^2(\F)=\bigoplus_{i+j=2}H^i(\F_a)\otimes H^j(\Ss^3)=0.\]

The condition $H^2(\F)=0$ can be interpreted as an obstruction to the
existence of a leafwise symplectic form on a foliation which exhibits
compactness. In order to explain this, let us recall that a Poisson
tensor $\pi\in\X^2(M)$ is said to be \emph{exact} if its Poisson
cohomology class $[\pi]\in H^2_\pi(M)$ is trivial. Note that this
means that:
\[ \Lie_X \pi=\pi,\]
for some vector field $X\in\X(M)$. Now we have:

\begin{proposition}
Let $\omega\in\Omega^2(\F)$ be a leafwise symplectic form. If the
class $[\omega]\in H^2(\F)$ vanishes then the associated Poisson
tensor $\pi$ is exact.
\end{proposition}

\begin{proof}
Let $\al\in\Omega^1(M)$ be a 1-form and denote by $\al^\sharp\in\X(M)$
the vector field obtained by contraction with the Poisson tensor
$\pi$. A simple computation shows that:
\[ \Lie_{\al^\sharp}\pi(\d f,\d g)=\d \al(\d f^\sharp,\d g^\sharp).\]
Therefore, if $\omega=\d_\F\al$, we conclude that:
\[ \Lie_{\al^\sharp}\pi(\d f,\d g)=\omega(\d f^\sharp,\d g^\sharp)
=\pi(\d f,\d g).\]
Hence, the vector field $X=\al^\sharp$ satisfies $\Lie_X \pi=\pi$,
and $\pi$ is exact.
\end{proof}

Recall that an \emph{integrable Poisson manifold} $(M,\pi)$ is
a Poisson manifold which arises as the space of units of a symplectic
groupoid $\G\tto M$ (see \cite{Wein2}).  The results in
\cite{CrFe1,CrFe2} show that the Poisson structures in this class are
the Poisson structures which have non-singular variations of
symplectic areas in directions transverse to the leaves. Simple
conditions on the foliation $\F$ will guarantee that \emph{any}
leafwise symplectic form on $\F$ determines an \emph{integrable}
Poisson tensor. For example, if $\F$ is a foliation such that:
\begin{enumerate}
 \item[(i)] $\F$ has no vanishing cycles;
 \item[(ii)] $\pi_2(L)$ is finite for every leaf $L$ of $\F$;
\end{enumerate}
then any Poisson structure with symplectic foliation $\F$ will be
integrable by an Hausdorff symplectic groupoid (see \cite{CrFe2}).

For an integrable Poisson manifold $(M,\pi)$, we will denote by
$\G(M)\tto M$ its source 1-connected symplectic groupoid, and we let
$\Omega\in\Omega^2(\G(M))$ denote its symplectic form. The following
proposition relates exactness of $\pi$ and of $\Omega$:

\begin{proposition}
\label{prop:exact}
An integrable Poisson structure $\pi$ on $M$ is exact if and only if
the symplectic form $\Omega$ on the associated symplectic groupoid is exact. 
\end{proposition}

\begin{proof}
Assume that $\pi$ is exact, so there exists a vector field $X\in\X(M)$
such that $\Lie_X\pi=\pi$. Then one checks easily that the flow
$\phi^t$ of $X$ satisfies: 
\[ (\phi^t)_*\pi=e^t\pi.\]
This means that, for each $t$, the map
\[ \phi^t:(M,\pi)\to (M,e^t\pi)\] 
is a Poisson diffeomorphism. By Lie's second theorem for Lie
algebroids (see \cite{CrFe1}), this integrates to a Lie groupoid
isomorphism:
\[
\xymatrix{
\G(M,\pi)\ar[r]^{\Phi^t}\ar@<1ex>[d]\ar@<-1ex>[d]& 
\G(M,e^t\pi)\ar@<1ex>[d]\ar@<-1ex>[d]&\\
M\ar[r]_{\phi^t}&M
}
\]
From the description of $\G$ given in \cite{CrFe1,CrFe2}, one sees
that the Lie groupoids $\G(M,\pi)$ and $\G(M,e^t\pi)$ can be
identified as Lie groupoids and they only differ on their symplectic
forms $\Omega$ and $\Omega_t$, which satisfy:
\[ \Omega_t=e^t\Omega.\]
In other words, we have that $\Phi^t:\G(M)\to \G(M)$ satisfies:
\[ (\Phi^{-t})^*\Omega=e^t\Omega.\]
If $\tilde{X}$ denotes the vector field on $\G(M)$ with flux
$\Phi^{-t}$, we conclude that:
\[ \Omega=\Lie_{\tilde{X}}\Omega=\d i_{\tilde{X}}\Omega,\]
so $\Omega$ is exact. The converse can be proved by reversing the
argument. 
\end{proof}

We have not used the fact that the symplectic form $\Omega$ on a
symplectic groupoid $\G$ is also \emph{multiplicative}: if
$m:\G^{(2)}\to\G$ is the multiplication in $\G$, and $s,t:\G\to M$ are
the source and target maps, then:
\[ m^*\Omega=s^*\Omega-t^*\Omega.\]
Using this condition one can show that (see \cite{Cr2}, Corollary 5.3):
\begin{proposition}
The symplectic form $\Omega$ is exact if and only if the restriction of
$\Omega$ to each $s$-fiber is exact.
\end{proposition}

Therefore, we see that that the condition $H^2(\F)=0$ places
strong restrictions on the topology of the groupoid integrating any
leafwise symplectic form on $\F$.

\begin{remark}
A different setting where one can discuss all these questions,
including Bertelson's notion of an open foliation, is the theory of
twisted Dirac brackets (see \cite{BCWZ}). In fact, given a foliation $\F$
of a manifold $M$, a leafwise 2-form $\omega\in\Omega^2(\F)$ defines a
$\phi$-twisted Dirac bracket on $M$, where $\phi=\d\omega$. Under some
nice assumptions, this twisted Dirac structure can be integrated to a
presymplectic groupoid $(\G,\Omega)$, where the presymplectic form $\Omega$
is $\phi$-twisted:
\[ \d\Omega=s^*\phi-t^*\phi.\]

When one varies the leafwise 2-form, the multiplication in the
groupoid and the presymplectic form changes but the total space and
the source and target maps stay the same. Hence, one can try to find a
homotopy, preserving the fibers, and which changes the presymplectic
form $\Omega$ to a symplectic form, and this should depend on 
topological properties of $\G$.
\end{remark}

All what we have said above applies only to \emph{regular} Poisson
structures. If one fixes a manifold $M$ and a \emph{singular}
foliation $\F$, one can look at the family of bivector fields
$\theta\in\X^2(M)$ which generate this foliation: each bivector
$\theta$ determines a bundle map $T^* M\to TM$ and we are interested
in those $\theta$ for which the image is $T\F$. Then one asks if,
given such an $\F$-compatible $\theta$, is it possible to find a
homotopy $\theta_t$ of $\F$-compatible bivector fields with
$\theta_0=\theta$ and $\theta_1$ a Poisson bivector field. Nothing
seems to be known about this more general problem.

\section{Rigidity}
\label{sec:rigid}
Let us now turn to \emph{rigidity} in Poisson geometry. 
Our main motivation is the well-known stability theorem in symplectic
geometry due to Moser:

\begin{theorem}
Let $M$ be a compact manifold, and let $\omega_t\in\Omega^2(M)$ be a
homotopy of symplectic structures in $M$ such that its cohomology
class $[\omega_t]\in H^2(M)$ is constant. Then $\omega_0$ and
$\omega_1$ are symplectomorphic.
\end{theorem}

A proof can be found in any standard text in symplectic geometry ,
such as the monograph of MacDuff and Salamon (\cite{MaSa}). We would
like to have analogues of this result in Poisson geometry. The first
difficulty we face is to find the appropriate compactness
assumption.

Recall that a Lie groupoid $\G\tto M$ is called \emph{proper} if the
map $(s,t):\G\to M\times M$ is a proper map. We propose the following
definitions:

\begin{definition}
\label{defn:proper}
A Poisson manifold $(M,\pi)$ is called \textbf{proper} if it is
integrable and its symplectic groupoid $\G(M)\tto M$ is proper. A
Poisson manifold $(M,\pi)$ is of \textbf{compact type} if it is 
compact and proper. 
\end{definition}

Note that $(M,\pi)$ is of compact type if and only it is
integrable and its source 1-connected groupoid $\G(M)$ is
compact. 

\begin{example}
A compact Poisson manifold may not be of compact type. For example, take any
compact Poisson manifold with the zero Poisson bracket: its symplectic
groupoid is $T^*M\to M$ (source and target coincide, and the product
is addition on the fibers), which is never proper.
\end{example}

\begin{example}
An example of a proper Poisson manifold is given by the linear Poisson
manifold $\gg^*$, where $\gg$ is a compact semi-simple Lie algebra:
its symplectic groupoid is $T^*G\tto\gg^*$, where $G$ is the (compact)
simply connected Lie group integrating $\gg$. Here the source/target
maps are given by right/left translations to the identity. This
Poisson manifold is not of compact type, since $\gg^*$ is not a
compact manifold.
\end{example}

\begin{example}
An example of a Poisson manifold of compact type is provided by a compact
symplectic manifold with finite fundamental group: its symplectic
groupoid is the fundamental groupoid $\pi_1(M)\tto M$, with
source/target maps given by the initial/end points of the homotopy
class of a path.
\end{example}

For a proper Poisson manifold we have the following properties:
\begin{proposition}
Let $(M,\pi)$ be a proper Poisson manifold. Then:
\begin{enumerate}
\item[(i)] All isotropy Lie algebras are of compact type.
\item[(ii)] All symplectic leafs are closed submanifolds.
\item[(iii)] Every Poisson vector field is Hamiltonian:
  $H^1_\pi(M)=0$.
\item[(iv)] There is a measure $\mu$ on $M$ which
is invariant under all Hamiltonian diffeomorphisms.
\end{enumerate}
\end{proposition}
\begin{proof}
Assume that $\G(M)\tto M$ is proper so that the map $(s,t):\G(M)\to
M\times M$ is a proper map. If $x\in M$ is a fixed point, its isotropy
Lie algebra $\gg_x$ integrates to the Lie group $G_x=(s,t)^{-1}(x,x)$
which is 1-connected and compact, so (i) holds.

Let $x_0\in M$ and let $L$ be the leaf through $x_0$. If $x_n\in L$ is
sequence converging to some $x\in M$, then there are elements
$g_n\in\G$ such that $(s,t)(g_n)=(x_0,x_n)$. The set
\[ K=\{(x_0,x_n):n\in\Nn\}\cup \{(x_0,x)\}\subset M\times M \]
is compact and $g_n\in (s,t)^{-1}(K)$. Therefore, the sequence $g_n$
has a convergent subsequence: $g_{n_i}\to g$. Observe that:
\begin{align*}
s(g)&=\lim s(g_{n_i})=x_0,\\
(x_0,x)&=\lim (s,t)(g_{n_i})=(s,t)(g).
\end{align*}
Hence $x\in L$ and $L$ is a closed leaf, so (ii) holds.

To prove (iii), we use the fact that (\cite{Cr}):
\[ H^1_\pi(M)=H^1_d(\G),\]
where $H^k_d(\G)$ denotes the groupoid cohomology with differentiable
cochains. Since $\G$ is proper, we have $H^1_d(\G)=0$, and (iii)
follows. In particular, the modular class must vanish and so (iv) also holds.
\end{proof}
 
Proper Poisson manifolds already exhibit a lot of rigidity, as it is
shown by the following result which is a consequence of general results
on linearization of proper groupoids due to Weinstein \cite{Wein3} and
Zung \cite{Zun1}:

\begin{theorem}
Let $(M,\pi)$ be a proper Poisson manifold, and let $L$ be a compact
symplectic leaf. Then $M$ is linearizable around $L$. 
\end{theorem}

We refer to \cite{Wein3,Zun1} for exact statements and proofs.
\vskip 10 pt

In order to obtain \emph{global} rigidity results one should require
some global compactness, and it is natural to consider Poisson
manifolds of compact type. At present, the only examples of
Poisson manifolds of compact type which we are aware of are
symplectic. Nevertheless, let us list some properties that any such
manifold must satisfy. 

We start with the following basic property, which shows that Poisson
manifolds of compact type extend compact symplectic manifolds:

\begin{proposition}
\label{prop:compt:exact}
If $(M,\pi)$ is a Poisson manifold of compact type then the Poisson
cohomology class $[\pi]\in H^2_\pi(M)$ is non-trivial. 
\end{proposition}

\begin{proof}
If $\pi$ was exact then, by Proposition \ref{prop:exact}, the
symplectic form $\Omega$ in $\G(M)$ would be exact, which contradicts
the fact that $\G(M)$ is compact. 
\end{proof}

\begin{remark}
Note that there are proper Poisson manifolds with $H^2_\pi(M)=0$. For
example, this happens for the Lie-Poisson structure on $M=\gg^*$, when 
$\gg$ is semi-simple of compact type.
\end{remark}

The following somewhat surprising property shows that the class of
Poisson manifolds of compact type is quite rigid:

\begin{proposition}
Let $(M,\pi)$ be a Poisson manifold of compact type. Then $\pi$
does not have any fixed points.
\end{proposition}

\begin{proof}
Assume that $\pi$ has a fixed point $x_0\in M$. Then 
$s^{-1}(x_0)=t^{-1}(x_0)$ is a compact, 1-connected Lie group. Hence,
it must be homologically 2-connected. By stability, all $s$-fibers are
homologically 2-connected. 

It was shown in \cite{Cr} that the classical Van Est homomorphism can
be extended to Lie groupoids: it relates the Lie algebroid
cohomology, the groupoid cohomology, and the ordinary cohomology of
the source fibers; when applied to Poisson manifolds it gives a map:
\[ \Phi: H^k_d(\G)\to H^k(T^*M)=H^k_\pi(M).\]
This map is an isomorphism up to degrees $k\le n$,
and injective in degree $n+1$, provided the source fibers are
homologically $n$-connected. 

In the present situation, we conclude that:
\[ H^2_\pi(M)=H^2_d(\G).\] 
Another basic result of \cite{Cr}, states that the differentiable groupoid
cohomology of a proper Lie groupoid vanishes (this is a
generalization of a well-known result in Lie theory stating that the
group cohomology of a compact Lie group vanishes). We
conclude that the class $[\pi]=0$, which contradicts Proposition
\ref{prop:compt:exact}. 
\end{proof}

{To} construct examples of Poisson manifolds of compact type one could
try to start with a proper Poisson manifold and restrict to a compact
Poisson submanifold. However, one must be careful since a Poisson
submanifold of an integrable Poisson manifold may be non-integrable,
as shown by the following example:

\begin{example}
Let $\gg$ be a semi-simple Lie algebra of compact type of dimension
$d+1$. We take $M=\gg^*$ with the Lie-Poisson bracket. The (negative
of the) Killing form defines an inner product on $\gg^*$ which is
$Ad^*$-invariant. It follows that the unit sphere
$\Ss^d=\{\xi\in\gg^*:||\xi||=1\}$ is a Poisson submanifold of
$\gg^*$. 

The Weinstein groupoid (see \cite{CrFe2}) of $\Ss^d$, which we still
denote by $\G(\Ss^d)$, can be obtained as follows: consider the
symplectic groupoid of $\gg^*$, namely $T^*G\tto\gg^*$, where $G$ is
the (compact) 1-connected Lie groupoid integrating $\gg$. Restricting 
to $\Ss^d$, we obtain the presymplectic groupoid $S(T^*G)\tto \Ss^d$,
where $S(T^*G)$ is the sphere bundle. This groupoid is still smooth,
but to obtain the Weinstein groupoid of $\Ss^d$ we have to factor out
the kernel of the presymplectic form. This may lead to a non-smooth
groupoid. For example, one can check that this happens if
$\gg=\mathfrak{su}(n)$, for $n\ge 3$. 
\end{example}

Another possible method to produce Poisson manifolds of compact type
is suggested by the following remark: Let $(S,\omega)$ be a symplectic
manifold and let $G\times S\to S$ be a proper and free action of a Lie
group $G$ by symplectomorphisms. Then $M=S/G$ is an integrable Poisson
manifold. In fact, its symplectic groupoid $\G(M)$ is obtained as the
symplectic quotient of the symplectic groupoid of $S$:
\[ \G(S/G)=\G(S)//G.\]
For a proof of this fact and generalizations to singular cases, we
refer to the upcoming paper \cite{FeOrRa}. If one starts with a
compact symplectic manifold $(S,\omega)$ with finite fundamental
group, this could lead to a Poisson manifold of compact type. However,
any canonical action on a compact symplectic manifold with finite
fundamental group cannot be free. So one is forced to consider actions
with a fixed isotropy type, and this makes this method much harder to
work.
\vskip 10 pt

In the regular case, we can state the following rigidity result:

\begin{proposition}
Let $\F$ be a regular foliation of a manifold $M$. Let
$\omega_t\in\Omega^2(\F)$ be a homotopy of foliated symplectic
structures such that its cohomology class $[\omega_t]\in H^2(\F)$ is
constant and the associated Poisson structures $\pi_t$ are of compact
type. Then $\pi_0$ and $\pi_1$ are Poisson diffeomorphic.
\end{proposition}

\begin{proof} We sketch a proof of this result.

Recall that there is a homomorphism $\d_\nu:H^2(\F)\to H^2(\F;\nu^*)$,
which can defined as follows: given a class $[\omega]\in H^2(\F)$
represented by a foliated 2-form $\omega\in\Omega^2(\F)$, pick an extension
$\widetilde{\omega}\in\Omega^2(M)$ of $\omega$. Since
$\d\widetilde{\omega}|_{T\F}=0$, we obtain a well-defined map
$\Gamma(\wedge^2T\F)\to\Gamma(\nu^*)$ given by:
\[ (X,Y)\mapsto \d\widetilde{\omega}(X,Y,\cdot).\]
One checks easily that this is a closed foliated 2-form with values in
$\nu^*$, whose cohomology class $\d_\nu[\omega]\in H^2(\F;\nu^*)$ does
not depend on any choices.

Now, if $\omega_t\in\Omega^2(\F)$ is a homotopy of foliated symplectic
structures such that the cohomology class $[\omega_t]\in H^2(\F)$ is
constant, the class $\d_\nu[\omega_t]\in H^2(\F,\nu^*)$ is independent
of $t$. The results of \cite{CrFe2} show that the groupoid
$\G(M,\pi_t)$ only depends on the class $\d_\nu[\omega_t]$, so we see
that the Poisson structures $\pi_t$ integrate to the same Lie groupoid
$\G(M)$, but with varying symplectic structures $\Omega_t$. Since
$\G(M)$ is compact, one can apply Moser's trick to produce an isotopy
of Lie groupoid automorphisms $\Phi_t$ such that
$(\Phi_t)^*\Omega_t=\Omega$. This isotopy covers the desired isotopy
of Poisson diffeomorphisms.
\end{proof}

Again in the non-regular case, just like in the previous section,
nothing seems to be known. In fact, it is not clear how such a general
result should be stated since there is no obvious way of comparing
Poisson cohomology classes attached to distinct Poisson structures. In
face of these difficulties, we now turn to milder rigidity
properties of Poisson structures.

\section{Stability of leaves}
\label{sec:stab}

The space $\Po(M)$ of Poisson structures on a manifold $M$ is a subset of the
space of smooth bivector fields $\X^2(M)$:
\[ \Po(M)=\set{\pi\in \X^2(M):[\pi,\pi]=0}.\]
Since $\X^2(M)$ is the space of sections of the vector bundle
$\wedge^2 TM$, we can consider on it the $C^k$ ($0\le k\le +\infty$)
compact-open topology. This induces topologies on $\Po(M)$ which, in
general, are extremely complicated. 

At the formal level, around a Poisson structure $\pi\in\Po(M)$, one
knows that the second Poisson cohomology group $H^2_\pi(M)$ is the
space of infinitesimal deformations of Poisson structures. In other
words, let $\M$ denote the moduli space of Poisson structures on $M$,
obtained by factoring out diffeomorphic Poisson structures:
\[ \M=\Po(M)/\Diff(M).\]
Then, formally, one has (see, e.g., \cite{CaWe}):
\[ T_\pi \M=H^2_\pi(M).\]
Hence, it is natural to study points $\pi$ where this tangent space
vanishes, where one hopes to have rigidity. In other words, one
is led to the question:
\begin{itemize}
\item If $H^2_\pi(M)=0$ is every nearby Poisson structure
  diffeomorphic to $\pi$? 
\end{itemize}
Of course, this question should be taken with care since
$H^2_\pi(M)$ only parametrizes the formal deformations of $\pi$. 

Let us take, for example, $M=\gg^*$ and assume that $\gg$ is of
compact type, so that $H^2_\pi(\gg^*)=0$. As a first step towards
understanding stability, we recall Conn's Linearization Theorem
\cite{Conn2}: 

\begin{theorem}
Let $\pi\in\Po(M)$ be a Poisson structure, $x_0\in M$ such that
$\pi_{x_0}=0$, and assume that the isotropy Lie algebra $\gg_{x_0}$
has compact type. Then there exists a neighborhood $V$ of $x_0$ which
is Poisson diffeomorphic to a neighborhood of zero in $\gg_{x_0}^*$.
\end{theorem}

Hence, if one could prove that every Poisson structure $\pi$ nearby
$\gg^*$ must have a zero $x_0$ with isotropy Lie algebra
$\gg_{x_0}\simeq\gg$, then, around $x_0$, $\pi$ would look like its
linearization $\gg^*$.  Therefore, one is led naturally to the study
of the \emph{stability} of zeros of Poisson structures:

\begin{definition}
A fixed point $x_0$ of a Poisson structure $\pi\in \Po(M)$ is
said to be \textbf{stable} if, given any neighborhood $V$ of $x_0$ in
$M$, there is a neighborhood $U$ of $\pi$ in $\Po(M)$, such that each
Poisson structure $\theta\in U$ has a fixed point in $V$.
\end{definition}

In other words, a fixed point of $\pi$ is stable if all nearby Poisson
structures have nearby fixed points. We have the following criterion
for stability of fixed points:

\begin{theorem}
\label{thm:stable:fxd:pt}
Let $\pi$ be a Poisson structure with a fixed point $x_0\in M$, and
denote by $\gg_{x_0}$ the isotropy Lie algebra at $x_0$. If
$H^2(\gg_{x_0})=0$ then $x_0$ is a stable fixed point.
\end{theorem}

For a proof of this result we refer to \cite{CrFe4}. In fact, there we
give a more precise version of the theorem showing that, under the
assumption $H^2(\gg_{x_0})=0$, the set of fixed points of $\pi$ in a
neighborhood of $x_0$ is a submanifold tangent to the fixed point set
of the linear approximation to $\pi$ at $x_0$, which is stable under
perturbations of the Poisson structure.

The condition in the theorem is, in some sense, not too far from
being a necessary condition. This can be seen by considering the case
of \emph{linear} Poisson structures. For these we can prove:

\begin{proposition}
\label{prop:stable:linear}
Let $\gg$ be a finite dimensional Lie algebra and take $M=\gg^*$ with
the linear Lie-Poisson bracket. If the origin is a stable fixed point
then $H^2(\gg)=0$.
\end{proposition}

\begin{proof}
Let us denote by $\pi_0$ the Lie-Poisson bracket on $\gg^*$, so that:
\[ \pi_0(\d f,\d g)(\xi)=-\langle \xi,[\d_\xi f,\d_\xi g]\rangle,\]
for any $\xi\in\gg^\ast$ and $f,g\in C^\infty(\gg^\ast)$. A cohomology
class $[c]\in H^2(\gg)$ is represented by a skew-symmetric map
$c:\gg\times\gg\to\Rr$. Consider the associated (constant) bivector
field on $\gg^\ast$ defined by:
\[ \pi_1(\d f,\d g)(\xi)=c(\d_\xi f,\d_\xi g),\]
for any $\xi\in\gg^\ast$ and $f,g\in C^\infty(\gg^\ast)$. Since
$\pi_1$ is constant, it defines a Poisson bivector field, i.e., the
Schouten bracket with itself vanishes: $[\pi_1,\pi_1]=0$.
Denote by $\delta:\wedge^\bullet\gg^*\to \wedge^{\bullet+1}\gg^*$ the
Chevalley-Eilenberg differential. Then:
\[ \delta c=0\quad \Longleftrightarrow\quad [\pi_0,\pi_1]=0,\]
and we conclude that the 1-parameter family of bivector fields
$\pi_t=\pi_0+t\pi_1$ is a pencil of Poisson structures.

Assume now that the origin is a stable fixed point point for $\pi_0$. 
Then for $t>0$, $\pi_t$ must also have a fixed point $\xi\in\gg^\ast$,
and we find:
\[
\pi_t(\xi)=\pi_0(\xi)+t\pi_1(\xi)=0\quad \Longleftrightarrow \quad 
tc=\delta\xi,\]
so $c$ must be a coboundary. Hence, if the origin is a stable fixed
point, we must have $H^2(\gg)=0$.
\end{proof}

As a corollary of Proposition \ref{prop:stable:linear} and Theorem
\ref{thm:stable:fxd:pt} we obtain:

\begin{corollary}
For a linear Poisson structure $M=\gg^*$ the origin is a stable fixed
point if and only if $H^2(\gg)=0$.
\end{corollary}

\begin{example}
\label{ex:semisimple}
Let $\pi$ be a Poisson structure on a manifold $M$, and let $x_0$ be a
fixed point of $\pi$ for which the isotropy Lie algebra $\gg_{x_0}$ is
semisimple. By the Second Whitehead Lemma, we have $H^2(\gg_{x_0})=0$,
so $x_0$ is a stable fixed point. 
\end{example}

\begin{remark} It was pointed out to us by Jean-Paul Dufour, that our
  criterion for fixed points can be generalized to higher degree
  singularities (see \cite{DuWa}).
\end{remark}

Fixed points of a Poisson structure are symplectic leaves of dimension
zero. The problem of stability of fixed points naturally generalizes
to higher dimensional symplectic leaves. This is made precise by the
following definition:

\begin{definition}
A symplectic leaf $S\subset M$ of a Poisson structure $\pi$ is said to be
a \textbf{stable leaf} if, for every tubular neighborhood $V$ of $S$ in
$M$, there is a neighborhood $U$ of $\pi$ in $\Po(M)$, such that each Poisson
structure $\theta\in U$ has a leaf contained in $V$ which is mapped
diffeomorphically onto $S$ by the tubular neighborhood projection.
\end{definition}

In the case where $\dim S=0$, i.e., $S=\{x_0\}$ is a fixed point, this
definition reduces to the definition of stability of a fixed point
given before. 

We conjecture that our criterion for stability of fixed points
(Theorem \ref{thm:stable:fxd:pt}) can be generalized to leaves. For
that recall the \emph{relative Poisson cohomology} of a symplectic
leaf $S$ of $(M,\pi)$ due to Ginzburg and Lu (see \cite{GiLu}): one
considers the complex $(\X^k(M,S),\delta)$ where:
\[ \X^k(M,S)=\Gamma(\wedge^k T^*_S M),\]
is the space of $k$-multivector fields along the symplectic leaf $S$,
and $\delta:\X^\bullet(M,S)\to\X^{\bullet+1}(M,S)$ is given by:
\[ \delta \theta=[\theta,\pi].\]
Since $\pi$ is tangent to $S$, one checks easily that $\delta$ is
well-defined. On the other hand, the vanishing of the Schouten bracket
$[\pi,\pi]=0$ implies that $\delta^2=0$, so $\delta$ is a
differential. One calls the cohomology of the resulting complex the
\textbf{relative Poisson cohomology} of the leaf $S$, and denotes it
by $H^\bullet_{\pi}(M,S)$. In terms of Lie algebroids,
$H^\bullet_{\pi}(M,S)$ coincides with the Lie algebroid cohomology of
the restriction $T^*_SM$.

\begin{example}
If $\dim S=0$, so that $S=\{x_0\}$ is a fixed point, it follows from
the definitions that the relative Poisson cohomology is just the
Lie algebra cohomology of the isotropy Lie algebra at $x_0$:
$H^\bullet_{\pi}(M,S)=H^\bullet_S(\gg_{x_0})$.
\end{example}

In this example, $H^\bullet_{\pi}(M,S)$ is finite dimensional. In fact,
one can show that the complex $(\X^k(M,S),\delta)$ is elliptic, hence
the relative Poisson cohomology is finite dimensional provided $S$ is of
finite type (e.g., if $S$ is compact). This is in sharp contrast with
the absolute Poisson cohomology, which is often infinite dimensional.

Now we pose the following natural conjecture:

\begin{conjecture}
Let $\pi$ be a Poisson structure with a compact symplectic leaf
$S\subset M$. If the 2nd Poisson cohomology relative to $S$ vanishes,
i.e., $H^2_\pi(M,S)=0$, then $S$ is a stable symplectic leaf.
\end{conjecture}

Further evidence for this conjecture is given in \cite{CrFe4}.
\vskip 10 pt

One may consider stronger notions of stability. For example, we may call
a symplectic leaf $S$ \textbf{strongly stable} if all nearby Poisson
structures have nearby symplectic leaves which are symplectomorphic
(rather than just diffeomorphic) to $L$, or one can even require the
restricted Lie algebroids $T^*_SM$ to be isomorphic. In the case of
fixed points, this would mean that nearby fixed points should have
isomorphic isotropy Lie algebras. These notions of stability seem
to be much harder to study.

\begin{remark}
For group (or Lie algebra) actions on a manifold $M$ the situation is
similar, as one may consider, also, two notions of stability: an orbit
$O\subset M$ of an action is called \emph{stable} (respectively,
\emph{strongly stable}) if every nearby action has a nearby leaf which
is diffeomorphic to (respectively, has the same orbit type as)
$S$. While the stability of orbits has been studied extensively
(cf.~\cite{St}), little seems to be known about strong stability.

Similarly, for foliations of a manifold $M$ of codimension $k$, one
also has two notions of stability: a leaf $L\subset M$ of a foliation
$\F$ is called \emph{stable} (respectively, \emph{strongly stable}) if
every nearby foliation has a nearby leaf which is diffeomorphic to
(respectively, has the same holonomy as) $L$. While the stability of
leafs has been studied extensively (cf.~the classical stability
results of Reeb \cite{Rb}, Thurston \cite{Th}, Langevin and Rosenberg
\cite{LaRo}), little seems to be known about strong stability.
\end{remark}

To finish, motivated by Conn's linearization result, we would like to
formulate a conjecture on linearization around leafs in the sense of
Vorobjev (\cite{Vor}):

\begin{conjecture}
Let $\pi$ be a Poisson structure with a compact symplectic leaf
$S\subset M$. If the restricted Lie algebroid $T^*_SM$ integrates to a
compact, source 1-connected, Lie algebroid, then $\pi$ is linearizable
around $S$.
\end{conjecture}

One should be able to prove this conjecture using the methods of
\cite{CrFe3}.  One can also look into milder properties of a leaf
which are shared by (some) nearby leaves. For example, given a leaf
$L$, one may ask if all nearby Poisson structures have a a nearby leaf
of the same dimension, etc. We refer the reader to the upcoming paper
\cite{CrFe4} for a detailed discussion of the problem of stability of
symplectic leaves.


\end{document}